\documentclass[a4paper, 11pt]{article}
\usepackage{amsfonts}
\usepackage{amsmath}
\usepackage{amssymb}
\usepackage{amsthm}

\usepackage[all]{nowidow}
\usepackage{logreq}
\usepackage{tikz}
\usepackage{paralist}	 
\usepackage{comment}
\usepackage{enumitem}
\usepackage{cite}
\usepackage{placeins}
\usepackage{multicol}
\usepackage{caption}
\usepackage{subfigure}
\usepackage[latin1]{inputenc}
\usepackage{graphicx}
\usepackage{epsf}
\usepackage{fancyhdr}
\usepackage{pstricks}
\usepackage{algorithm}
\usepackage{algorithmic}
\usepackage{lscape}
\newpsobject{showgrid}{psgrid}{subgriddiv=1, griddots=10,gridlabels=6pt}
\usepackage{authblk}
\usepackage{csquotes}
\usepackage[none]{hyphenat}
\usepackage[capitalise]{cleveref}
\usepackage[includefoot,bottom=15pt]{geometry}

\begin{document}

\begin{center} 
{A heuristic for the non-unicost set covering problem using local branching} ~\\ 
~\\
{J.E. Beasley} ~\\
~\\
~\\
Mathematics, Brunel University, Uxbridge UB8 3PH, UK 
 ~\\
 ~\\
 john.beasley@brunel.ac.uk 
~\\

{http://people.brunel.ac.uk/$\sim$mastjjb/jeb/jeb.html}
~\\
~\\
April 2023, Revised November 2023\\
\end{center}

\begin{abstract}
In this paper we present a heuristic for the non-unicost set covering problem using local branching.
Local branching   eliminates the need to define a problem specific search neighbourhood  for any particular (zero-one) optimisation problem. It does this by incorporating a generalised Hamming distance neighbourhood into the problem, and this leads naturally to an appropriate neighbourhood search procedure. 
We apply our approach  to the non-unicost set covering problem. Computational results are presented for 65 test problems that have been widely considered in the literature. Our results indicate that
our heuristic is better than six of the eight other heuristics we examined,  slightly worse than that of one heuristic, but that
there is a single heuristic that out-performs all others.
We believe that the work described here illustrates that the potential for using local branching, operating as a stand-alone matheuristic, 
has not been fully exploited in the literature.

\end{abstract}

\sloppy Keywords: Hamming distance;   local branching; integer programming; matheuristic; neighbourhood search; set covering

\section{Introduction}

As the reader may be aware a common approach to the heuristic solution of many zero-one integer programming problems is to apply neighbourhood search. By this we mean that given a (typically feasible) solution to the problem at hand we examine \enquote{small} changes to this solution. So we examine solutions in the \enquote{neighbourhood} of this feasible solution. If we find a better feasible solution then this typically becomes the new  solution and the process repeats until some termination condition is satisfied 
(e.g.~computational time limit, or failure to improve on the solution).

There are a number of general neighbourhood search approaches in the literature such as simulated annealing~\cite{kirkpatrick83},
tabu search~\cite{glover90} and 
variable neighbourhood search~\cite{mladenovic97,hansen01} that can be applied. Such approaches set out a general search procedure, but need particularisation for the problem at hand, e.g.~in defining the neighbourhood of a solution. Typically the neighbourhood of a solution is defined by specifying the possible moves away from a solution. These neighbourhood search approaches have been extensively used in the literature. For example a recent search using Web of Science (http://www.webofscience.com) listed approximately 33,000 papers referring to simulated annealing, 12,000 papers referring to tabu search and 7,000 papers referring to variable neighbourhood search

In this paper we present an optimisation based approach to neighbourhood search making use of local branching. 
Local branching eliminates the need to define a problem specific search neighbourhood for any particular (zero-one) optimisation problem. 
It does this by incorporating a generalised Hamming distance neighbourhood into the problem, and this leads naturally to an appropriate 
neighbourhood search procedure. 

The structure of this paper is as follows. In Section~\ref{opt} we define what we mean by the neighbourhood of a feasible solution to a general (zero-one) optimisation problem. We then go on to outline a search procedure that we can adopt to successively search for improved solutions.
 In Section~\ref{app} we consider the example optimisation problem, the non-unicost set covering problem, to which we are going to apply
 our
approach. We define the problem and consider relevant literature on the problem with especial reference to papers in the literature which report good computational results.
In Section~\ref{results} we present computational results produced by applying our 
approach to 65 non-unicost set covering problems  that have been extensively considered by others in the literature. 
 We compare our results with those obtained using Cplex~\cite{cplex1210} alone.
We also give a comparison between the results from our approach and eight other approaches 
presented previously in the literature.
Finally in Section~\ref{conc} we present our conclusions.

\section{Optimisation based neighbourhood search}
\label{opt}
In this section we first define 
 what we mean by the neighbourhood of a feasible solution to a general (zero-one) optimisation problem. 
We then go  on to outline a search procedure that we can adopt to successively search for improved solutions.

\subsection{Neighbourhood}
To illustrate our approach suppose that we have a general zero-one integer programming problem
involving $n$ zero-one variables $[x_i,~i=1,\ldots,n]$ and $m$ constraints where the optimisation problem is:

\begin{equation}
\mbox{minimise}~~~~\sum_{i=1}^n c_i x_i
\label{eq1}
\end{equation}
subject to:
\begin{equation}
\sum_{j=1}^n  a_{ij} x_j \geq b_i ~~~~ i = 1, \ldots, m
\label{eq2}
\end{equation}
\begin{equation}
x_i \in \{0,1\}  ~~~~i=1,\ldots,n
\label{eq3}
\end{equation}
Equation~(\ref{eq1}) is a minimisation objective. Without significant loss of generality we shall henceforth assume that all the objective function coefficients $[c_i]$ are integer. Equation~(\ref{eq2}) represents the constraints of the problem and Equation~(\ref{eq3}) the integrality condition.

Let $[X_i]$ be some feasible solution to the problem. Then in this paper we define the neighbourhood of  $[X_i]$ to be any set of zero-one variable values $[x_i]$ satisfying $1 \leq \sum_{i=1}^n |x_i - X_i| \leq K$, where $K$ is a known positive constant of our choice. In other words the neighbourhood of $[X_i]$ is any set of zero-one values $[x_i]$ such that the Hamming distance between $[x_i]$ and $[X_i]$ lies between one and $K$.

Then consider the optimisation problem:
\begin{equation}
\mbox{minimise}~~~~\sum_{i=1}^n c_i x_i
\label{eq1a}
\end{equation}
subject to Equations~(\ref{eq2}),(\ref{eq3}) and:
\begin{equation}
1 \leq \sum_{i=1~X_i=0}^n  x_i  + \sum_{i=1~X_i=1}^n (1- x_i ) \leq K
\label{eqs1}
\end{equation}
\begin{equation}
\sum_{i=1}^n c_i x_i \leq \sum_{i=1}^n c_i X_i - 1
\label{eqs2}
\end{equation}
Here we have added two constraints to our original optimisation problem. In Equation~(\ref{eqs1}) the expression seen is a linearisation of the nonlinear Hamming distance $\sum_{i=1}^n |x_i - X_i|$. 
Equation~(\ref{eqs1}) ensures that the Hamming distance between $[x_i]$ and $[X_i]$ is at least one (so we have a  solution different from $[X_i]$) and is also less than or equal to $K$.

 Equation~(\ref{eqs2}) implies that we are only interested in improved feasible solutions in the neighbourhood of $[X_i]$, i.e.~those that strictly improve on the solution value $\sum_{i=1}^n c_i X_i $ associated with the current solution. Improving on the current feasible solution cannot be guaranteed by Equation~(\ref{eqs1}) since it only constrains the structural (Hamming distance) difference between two solutions, it does not address their objective function values. 

With regard to a minor technical issue here we have that in integer programming terms use of Equation~(\ref{eqs2}) automatically implies that  the Hamming distance between $[x_i]$ and $[X_i]$ is at least one.
This is because any improved feasible solution must be different from $[X_i]$. However it could be that including an explicit lower limit on the Hamming distance of one (as in Equation~(\ref{eqs1})) improves computational performance, e.g.~by improving the linear programming relaxation solution, so we include it here.

\sloppy Our  amended optimisation problem is now
optimise Equation~(\ref{eq1a}) subject to Equations~(\ref{eq2}),(\ref{eq3}),(\ref{eqs1}),(\ref{eqs2}). In essence
here we have amended the original optimisation problem to restrict attention to distinctly different (and improved) solutions within the $K$ neighbourhood of $[X_i]$. 
Note here that because of the extra constraints added to the original optimisation problem  
any feasible solution to the amended optimisation problem 
must be an improved solution as compared to $[X_i]$.

Use of a constraint involving Hamming distance has previously been given in the literature by Fischetti and Lodi~\cite{fischetti03}. In their approach, which they call \enquote{local branching}, once a feasible solution $[X_i]$ is found within an enumerative scheme, for example linear programming based tree search, two tree branches are created. One of these branches, which they call the left branch, has 
$ \sum_{i=1}^n |x_i - X_i| \leq K$.
The other branch, which they call the right branch, has 
$ \sum_{i=1}^n |x_i - X_i| \geq K +1$.
They suggested tactical exploration of the left branch, using standard branching procedures, in the hope of finding an improved feasible solution within the Hamming distance $K$ neighbourhood of the current feasible solution before proceeding with exploration of the right branch.
They proposed varying the value of $K$ depending upon search progress: for example reducing $K$ if the left branch has not resulted in an improved solution within a specified time limit; increasing $K$ to diversify the search.

Our approach differs from local branching as described in~\cite{fischetti03} in one significant respect, namely that
we only focus on the left branch, no exploration is attempted with regard to the right branch. This is because we are focusing on generating good quality heuristic solutions, abandoning any attempt to achieve a provably optimal solution for the original problem (Equations~(\ref{eq1})-(\ref{eq3})) under investigation.

In general terms it is clear that if we solve our amended optimisation problem to proven global optimality (e.g.~using a package such as Cplex~\cite{cplex1210}) then we will either:
\begin{compactitem}
\item find an improved feasible solution within the $K$ neighbourhood of $[X_i]$, technically the minimum feasible solution within the neighbourhood; or
\item prove that there is no improved feasible solution within the neighbourhood.
\end{compactitem}
Obviously computational considerations may mean that we do not solve the amended optimisation problem to proven global optimality, but within computational limits we may still find an improved feasible solution. 
As noted above any feasible solution to the amended optimisation problem must by definition be an improved solution as compared to $[X_i]$.
 Obviously, as with standard neighbourhood search procedures, an improved feasible solution can be used to replace $[X_i]$ and the process repeated in a natural way. The search procedure we adopted making use of the amended optimisation problem is detailed below.

\subsection{Search procedure}

Our search procedure  requires an initial feasible solution $[X_i]$ as well as three parameter values. These are:
an initial value for $K$; a value $\delta$ for incrementing $K$ so as to increase the size of the neighbourhood and a value $L$ for the number of successive iterations we allow without improving the solution before terminating the search.
\newline
\newline
\noindent Our search procedure is:
\begin{enumerate}[label=(\alph*), noitemsep]
\item Initialise $[X_i]$, $K$,  $\delta$ and $L$. Set $t \leftarrow 0$, where $t$ is the iteration counter.
\item Set $t \leftarrow t+1$. Solve the amended optimisation problem. If we find an improved feasible solution replace   $[X_i]$ with this solution. 
\item If $L$ successive iterations have been performed without improving the current feasible solution then stop, else set $K \leftarrow K + \delta$  and go to step (b).
\end{enumerate}

\noindent In this procedure we increase the size of the neighbourhood (increment $K$ by $\delta$) at each iteration, irrespective as to  whether an improved solution has been found or not.  This ensures that we continually expand the search space around the (current) feasible solution.  The procedure only terminates once $L$ successive iterations have been performed without finding an improved feasible solution.

\subsection{Comment}
We would make a number of comments as to our 
approach: 
\begin{itemize}
\item Our approach  draws directly on the mathematical formulation of the problem and hence can be classed as a matheuristic~\cite{boschetti2009, boschetti23, maniezzo21}.
\item It  eliminates the need to design problem specific search neighbourhoods, since the neighbourhood is automatically incorporated into 
the amended optimisation problem using the Hamming distance (Equation~(\ref{eqs1})) as discussed  above. 
\item If the amended optimisation problem associated with the final feasible solution found has been solved to proven global optimality then we have an \emph{\textbf{absolute guarantee}} that there is no improved feasible solution within the Hamming distance $K$ neighbourhood associated with that final feasible solution.
\item Clearly local branching is not a new concept. Indeed various authors in the literature have mentioned its use as a heuristic, e.g.~most recently~\cite{boschetti23}.  \textbf{\emph{However,  we believe that work described here illustrates that the potential for using local branching, operating as a stand-alone matheuristic as in the approach described in this paper,
has not been fully exploited in the literature.}} 
\end{itemize}

\section{The set covering problem}
\label{app}

The set covering problem is the problem of choosing a minimum cost set of columns $[x_i]$ that collectively cover each of the $m$ rows in the problem. Referring back to Equation~(\ref{eq2}) above we have that $[a_{ij}]$ is a known matrix with $a_{ij}=1$ if column $j$ covers row $i$, $a_{ij}=0$ otherwise. The values $[b_i]$ are all one. 

There are two variants of the problem, one where the  column costs $[c_i]$ are all one (known as the unicost set covering problem), one where the column costs $[c_i]$ are general 
non-negative values (referred to as the non-unicost set covering problem, or more commonly as just the set covering problem). In the results given below we apply our approach to the 
non-unicost problem. Of the two variants of the problem the non-unicost variant has attracted greater attention in the literature.

The (non-unicost) set covering problem has been considered by a number of authors in the literature as discussed below.
It is not our intention here to give a comprehensive and detailed review of the literature for the set covering problem. Indeed that would be a mammoth task, since a recent search using Web of Science (http://www.webofscience.com) listed over 500 papers that included the phrase \enquote{set covering} in their title.
Rather our intention is to highlight significant papers in the literature which report good computational results on set covering instances. This is because the focus of this paper is whether, for the specific example problem (set covering) considered, our 
approach can yield good quality results as compared with those already reported in the literature. 

\subsection{Relevant literature}

\sloppy Caprara et al~\cite{caprara00} give a survey of algorithms for the set covering problem prior to 2000. As is clear from their paper most of the authors in the literature since 1990 have made use of the test problems publicly available from 
OR-Library~\cite{beasley1990}, see 
 http://people.brunel.ac.uk/$\sim$mastjjb/jeb/info.html. 
In this paper we also make use of these test problems.

Lan et al~\cite{lan07}, based on Lan~\cite{lan04}, presented a heuristic approach which they called Meta-RaPS (Meta-heuristic for Randomized Priority Search). They stressed the use of randomness to avoid local optima. Their approach is a repeated application of: firstly a constructive heuristic to find a feasible solution (but including randomisation); secondly a local improvement heuristic using neighbourhood search. Their approach also included preprocessing, both to exclude columns from consideration and to include them if a column is the only one that covers a row. To reduce the computation time associated with their neighbourhood search procedure they defined a core problem consisting of a small subset of columns. 

Lan et al~\cite{lan07} reported that their heuristic is one of only two to find all optimal/best-known solutions for non-unicost instances. 
They gave a table illustrating the effectiveness of different heuristics on 65 non-unicost set covering problems. Of note there is 
the indirect genetic algorithm of Aickelin~\cite{aickelin02};
 the genetic algorithm of Beasley and Chu~\cite{beasley96}; and
the lagrangian heuristic of Caprara et al~\cite{caprara99}.
We consider each of these three approaches below.

Aickelin~\cite{aickelin02} presented a genetic algorithm approach with a decoder which works on a permuted list of the rows to be covered, with hill-climbing to improve the solution applied after the decoder has provided a suitable solution. For 65 non-unicost set covering problems they compared their approach with other approaches,~\cite{beasley96,caprara99},  taken from the literature.

 Beasley and Chu~\cite{beasley96} presented a genetic algorithm approach including a new fitness-based crossover operator (fusion), a variable mutation
rate and a heuristic feasibility operator tailored specifically for the non-unicost set covering problem. They reported computational experience for
 their approach  on 65 non-unicost set covering problems.

Caprara et al~\cite{caprara99} presented a lagrangian heuristic approach using dynamic pricing for the variables and systematic use of column fixing to improve the solution. They made use of a number of improvements in the subgradient optimisation procedure as well as a refining procedure to improve upon any given solution.  They gave a table illustrating the effectiveness of different heuristics on 65 non-unicost set covering problems showing that their approach performs well.

Naji-Azimi et al~\cite{naji10} presented an electromagnetic metaheuristic approach drawing on the work of Birbil and Fang~\cite{birbil03}. Their approach involves an initial  preprocessing step and then repetitively adjusting a pool of solutions to which local search is first applied and where the solutions are then changed based on  the force generated by the  \enquote{charge} associated with each solution. Mutation was also applied to perturb solutions. They considered 65 non-unicost problems and compared their results with those of Lan et al~\cite{lan07}. Their results indicated that their approach was competitive with that of Lan et al~\cite{lan07}.

Reyes and Araya~\cite{reyes21} presented a greedy randomised adaptive search procedure (GRASP~\cite{feo95, festa02}) based strategy for the non-unicost set covering problem. They proposed iterated local search and reward/penalty procedures to accelerate convergence and  improve upon the GRASP solutions. Their approach also included preprocessing both to exclude columns from consideration and to include them. They presented results, based upon 30 trials,  for 65 non-unicost set covering problems.

In recent years a number of papers in the literature have applied algorithms based upon metaphors/paradigms drawn from the natural world (sometimes referred to as bio-inspired metaheuristics). One  example of work of this kind is Soto et al~\cite{soto17}  who presented approaches based on cuckoo search and black hole optimisation. However note here that metaheuristic work based on metaphors from nature  has its critics, e.g.~\cite{aranha22}.

Cuckoo search (see Yang and Deb~\cite{yang09}) is a population based approach  where each \enquote{nest} in the population contains a number of \enquote{eggs} (solutions) and \enquote{cuckoos} lay eggs in randomly chosen nests. The best nests carry over to the next generation.
Black hole optimisation (see Kumar et al~\cite{kumar15}) is also a population based approach where a \enquote{black hole} attracts \enquote{stars} (solutions). Stars change locations as they are attracted by the black hole. Some stars are absorbed by the black hole and replaced by newly generated stars (solutions).

Soto et al~\cite{soto17}  applied these two approaches to the non-unicost set covering problem. They applied preprocessing and presented results for 65 non-unicost set covering problems (based on 30 trials for each instance).

\section{Computational results}
\label{results}

In this section we first discuss the non-unicost set covering test problems which we used. We then give computational results for our 
approach when applied to these test problems. We compare our results with those obtained using Cplex~\cite{cplex1210} alone.
We also give a comparison between the results from our approach and eight other approaches 
presented previously in the literature.

\subsection{Test problems}

We used the standard set of 65 non-unicost set covering test problems that are  available from 
OR-Library~\cite{beasley1990}, see http://people.brunel.ac.uk/$\sim$mastjjb/jeb/info.html. 
  We used a Windows pc with 8GB of memory, a 512GB SSD disk and an Intel Core i5-1135G7 2.4GHz processor, a multi-core pc with four cores. The initial value of $K$ was set to 5, the neighbourhood increment $\delta$ was set to 5 and the number successive iterations $L$  allowed without improving the current feasible solution  before termination was set to 5. These values were set based on limited computational experimentation.

The initial feasible solution required to start the search was produced using a simple greedy heuristic for the set covering problem: repetitively choosing a column with the minimum value of the ratio (column cost/number of uncovered rows covered by the column). Once a solution covering all rows had been found we removed any redundant columns (those columns for which all rows which they cover are also covered by other columns).
Note here that, unlike a number of other papers in the literature 
(e.g.~\cite{lan04, lan07,naji10,reyes21,soto17}), we made no use of problem-specific preprocessing
to eliminate columns. 

Table~\ref{table1} shows the characteristics of the 65 non-unicost problems considered. In that table we show the problem set name, the number of instances in that set, the number of rows ($m$)  and columns ($n$) and the problem density ($[\sum_{i=1}^m \sum_{j=1}^n a_{ij}/mn]$ expressed as a percentage).  We set a time limit of 15 seconds for the solution of each optimisation 
 for all the problems with $m \leq 500$, where Cplex~\cite{cplex1210} with default parameter settings was used as the solver. 
 For the larger problems with $m=1000$ we increased this time limit to 45 seconds.  Each of the 65 test problems was only executed one, i.e.~we did not do a number of 
different replications, each time with a differing initial feasible solution.

\begin{table}[hbt!]
\centering
\renewcommand{\tabcolsep}{1mm} 
\renewcommand{\arraystretch}{1.0} 
\begin{tabular}{|c|c|c|c|c|}
\hline 
Problem set & Number of & Number of  & Number of 
 & Density\\
name & instances & rows ($m$) & columns ($n$) &  (\%) \\
\hline
4 & 10 & 200 & 1000 & 2 \\
5 & 10 & 200 & 2000 & 2 \\
6 & 5 & 200 & 1000 & 5 \\
A & 5 & 300 & 3000 & 2 \\ 
B & 5 & 300 & 3000 & 5 \\
C & 5 & 400 & 4000 & 2 \\ 
D & 5 & 400 & 4000 & 5 \\
NRE & 5 & 500 & 5000 & 10 \\ 
NRF & 5 & 500 & 5000 & 20 \\
NRG & 5 & 1000 & 10000 & 2 \\ 
NRH & 5 & 1000 & 10000 & 5 \\
\hline
\end{tabular}
\caption{Test problem characteristics}
\label{table1}
\end{table}

\subsection{Results}

Table~\ref{table2} shows the results obtained by our optimisation approach. 
In that table we show the optimal/best-known solution value (OBK) for each problem, as taken from Lan et al~\cite{lan07}. We also show the solution value as obtained by our approach, 
 the final value of $K$ at termination, 
the total computation time (in seconds),  and  whether the amended optimisation problem for the final value of $K$ was solved to proven optimality or not. We have not given in Table~\ref{table2} the number of  iterations made in our procedure as this can easily deduced by dividing the final value of $K$ by the neighbourhood increment $\delta$, where we used $\delta=5$.

So for example consider problem 4.10 in Table~\ref{table2}. The optimal/best-known solution for this instance is 514 and the \enquote{o} signifies that this was the solution found by our approach. The value of $K$, the size of the neighbourhood at the final iteration, was 45 and the total solution time was 0.4 seconds. The \enquote{yes} in the solution guarantee column signifies that the amended optimisation problem associated with this final value of $K$ was solved to proven optimality, indicating that we have an absolute guarantee that there is no improved solution within a Hamming distance of $K=45$ from the solution associated with the value of 514 as found by our approach.

Considering Table~\ref{table2} it is clear that for all but one of the 65 test problems we found the optimal/best-known solution. From Table~\ref{table1} these problems are of increasing size and for all 45 problems up to and including problem set D we have the guarantee on solution quality. However for the 20 larger problems we only have one instance in which this is the case (recall here that we impose a time limit for the solution of each and every amended optimisation problem encountered during the process).

\begin{table}[hbt!]
\footnotesize
\centering
\renewcommand{\tabcolsep}{1mm} 
\renewcommand{\arraystretch}{0.85} 
\begin{tabular}{|c|c|c|c|c|c|}
\hline 
Instance & Optimal/best-known (OBK) & Solution & Final $K$ & Total time (secs)
 & Solution guarantee \\
\hline
4.1	&	429	&	o	&	45	&	0.4	&	yes	\\
4.2	&	512	&	o	&	50	&	0.6	&	yes	\\
4.3	&	516	&	o	&	55	&	0.5	&	yes	\\
4.4	&	494	&	o	&	50	&	0.5	&	yes	\\
4.5	&	512	&	o	&	45	&	0.4	&	yes	\\
4.6	&	560	&	o	&	55	&	0.7	&	yes	\\
4.7	&	430	&	o	&	50	&	0.4	&	yes	\\
4.8	&	492	&	o	&	45	&	0.8	&	yes	\\
4.9	&	641	&	o	&	55	&	1.0	&	yes	\\
4.10	&	514	&	o	&	45	&	0.4	&	yes	\\
\hline
5.1	&	253	&	o	&	50	&	1.1	&	yes	\\
5.2	&	302	&	o	&	50	&	1.3	&	yes	\\
5.3	&	226	&	o	&	50	&	0.7	&	yes	\\
5.4	&	242	&	o	&	45	&	1.1	&	yes	\\
5.5	&	211	&	o	&	45	&	0.7	&	yes	\\
5.6	&	213	&	o	&	45	&	0.8	&	yes	\\
5.7	&	293	&	o	&	50	&	1.1	&	yes	\\
5.8	&	288	&	o	&	45	&	1.0	&	yes	\\
5.9	&	279	&	o	&	50	&	0.8	&	yes	\\
5.10	&	265	&	o	&	50	&	0.7	&	yes	\\
\hline
6.1	&	138	&	o	&	60	&	2.0	&	yes	\\
6.2	&	146	&	o	&	40	&	1.1	&	yes	\\
6.3	&	145	&	o	&	45	&	1.5	&	yes	\\
6.4	&	131	&	o	&	45	&	1.4	&	yes	\\
6.5	&	161	&	o	&	50	&	1.8	&	yes	\\
\hline
A1	&	253	&	o	&	50	&	3.7	&	yes	\\
A2	&	252	&	o	&	55	&	4.1	&	yes	\\
A3	&	232	&	o	&	50	&	2.9	&	yes	\\
A4	&	234	&	o	&	50	&	2.4	&	yes	\\
A5	&	236	&	o	&	65	&	3.3	&	yes	\\
\hline
B1	&	69	&	o	&	40	&	3.5	&	yes	\\
B2	&	76	&	o	&	45	&	9.4	&	yes	\\
B3	&	80	&	o	&	45	&	4.7	&	yes	\\
B4	&	79	&	o	&	45	&	10.2	&	yes	\\
B5	&	72	&	o	&	40	&	3.4	&	yes	\\
\hline
C1	&	227	&	o	&	55	&	4.6	&	yes	\\
C2	&	219	&	o	&	60	&	5.8	&	yes	\\
C3	&	243	&	o	&	80	&	19.6	&	yes	\\
C4	&	219	&	o	&	50	&	4.0	&	yes	\\
C5	&	215	&	o	&	55	&	5.0	&	yes	\\
\hline
D1	&	60	&	o	&	50	&	7.6	&	yes	\\
D2	&	66	&	o	&	50	&	23.0	&	yes	\\
D3	&	72	&	o	&	55	&	23.5	&	yes	\\
D4	&	62	&	o	&	45	&	10.6	&	yes	\\
D5	&	61	&	o	&	45	&	5.1	&	yes	\\
\hline
NRE1	&	29	&	o	&	40	&	78.5	&	no	\\
NRE2	&	30	&	o	&	65	&	160.2	&	no	\\
NRE3	&	27	&	o	&	60	&	139.8	&	no	\\
NRE4	&	28	&	o	&	50	&	105.3	&	no	\\
NRE5	&	28	&	o	&	40	&	76.8	&	no	\\
\hline
NRF1	&	14	&	o	&	45	&	114.2	&	no	\\
NRF2	&	15	&	o	&	35	&	78.5	&	no	\\
NRF3	&	14	&	o	&	45	&	44.4	&	yes	\\
NRF4	&	14	&	o	&	40	&	93.4	&	no	\\
NRF5	&	13	&	o	&	45	&	109.4	&	no	\\
\hline
NRG1	&	176	&	o	&	65	&	308.3	&	no	\\
NRG2	&	154	&	o	&	60	&	219.7	&	no	\\
NRG3	&	166	&	o	&	85	&	565.3	&	no	\\
NRG4	&	168	&	o	&	115	&	821.2	&	no	\\
NRG5	&	168	&	o	&	75	&	475.4	&	no	\\
\hline
NRH1	&	63	&	64	&	55	&	344.9	&	no	\\
NRH2	&	63	&	o	&	70	&	501.0	&	no	\\
NRH3	&	59	&	o	&	100	&	773.0	&	no	\\
NRH4	&	58	&	o	&	85	&	618.6	&	no	\\
NRH5	&	55	&	o	&	55	&	342.0	&	no	\\

\hline
\end{tabular}
\caption{Computational results}
\label{table2}
\end{table}
\normalsize

\subsection{Comparison with Cplex}

In order to compare our results with those obtained using Cplex alone we solved all of the 65 test problems, 
using Cplex~\cite{cplex1210} with default parameter settings, but imposed a time limit for each problem. We conducted two experiments here: firstly, with the time limit (TL)  for each problem set equal to the corresponding total solution time as shown in Table~\ref{table2}; secondly 
with the time limit (TL) for each problem set equal to 1000 seconds (so set to a time greater than the largest solution time seen in Table~\ref{table2}).
The results are shown in Table~\ref{table4}. In that table we give for each problem the OBK solution (as in Table~\ref{table2}). For each of the two Cplex experiments we give the solution value, where again an \enquote{o} signifies that solution found was equal to the OBK solution, together with the total computation time.

Over the 65 test problems where the time limit for Cplex was set equal to the total time taken by our heuristic as in Table~\ref{table2} 
we have that for 60 of the 65 test problems the solution found was equal to the OBK solution, corresponding to an average percentage deviation of 0.12\%, in an average total time of 79.7 seconds. 
Over the 65 test problems where the time limit for Cplex was set equal to 1000 seconds
we have that for 62 of the 65 test problems the solution found was equal to the OBK solution, corresponding to an average percentage deviation of 0.06\%, in an average total time of 146.4 seconds. 
The comparative values for the approach given in this paper as in Table~\ref{table2} are that 
for 64 of the 65 test problems the solution found was equal to the OBK solution, corresponding to an average percentage deviation of 0.02\%, in an average total time of 94.6 seconds. All of these results are summarised in the \enquote{All 65 problems} section of Table~\ref{table5}.

One point to note from Table~\ref{table4} is that for the results in the last two columns of Table~\ref{table4}, with a time limit of 1000 seconds, any problem solved within that time has been solved to proven optimality by Cplex. So here all problems up to and including problem NRF5 have been solved to proven optimality (with the maximum solution time for any problem being 41.0 seconds). By contrast for the last ten problems in Table~\ref{table4} only one problem is solved to proven optimality within 1000 seconds (namely NRG2 requiring 335.1 seconds).
This  is a result of the advances in both hardware and software (as well as optimisation theory)
since these set covering problem instances were first put forward in the literature during the period 
1980-1990~\cite{balas80, beasley87, beasley90} and made freely available in 1990~\cite{beasley1990}, so over 30 years ago. \emph{\textbf{However we have included all of 
these test problems in Table~\ref{table2} and Table~\ref{table4} since all other workers have also made use of these test problems in their published work, as discussed later below.}} 

Of course we might well take the view that there is little point in adopting an heuristic for any particular set covering problem that can be easily solved to proven optimality by a modern optimisation package such as Cplex. Rather the worth of a heuristic is given by how well it performs on problems that cannot so easily solved. To illustrate this we took just the nine NRG and NRH problems not solved to proven optimality in 1000 seconds and show in 
the associated section of Table~\ref{table5} the comparative values for the heuristic presented in this paper and Cplex alone. Note here that the averages shown in the two sections of  Table~\ref{table5} are averages over 65 problems in the first section, averages over 9 problems in the second section (which is why, for example, average percentage deviation is higher in the second section).

Considering Table~\ref{table5} it seems clear that our heuristic is adding value, resulting in a lower average percentage deviation and a higher number of optimal/best-known solutions as compared with using Cplex alone, particularly for the nine larger problems that are more challenging to solve using Cplex.

As far as we are aware no proven optimal values for the nine NRG and NRH problems not solved within a 1000 second time limit have been given in the literature. In order to resolve this issue, and provide known optimal values for future workers, we allowed Cplex unlimited time to solve these problems.
 The results are given in Table~\ref{tableopt}. On the pc we used we could solve all four NRG problems to 
proven optimality, however we had insufficient memory/disk space to solve any of the NRH problems to proven optimality.
As can be seen in Table~\ref{tableopt} all of these NRG problems required substantially more than 1000 seconds to solve to proven optimality.

\begin{table}[hbt!]
\footnotesize
\centering
\renewcommand{\tabcolsep}{1mm} 
\renewcommand{\arraystretch}{0.85} 
\begin{tabular}{|c|c|c|c|c|c|}
\hline 
Instance & Optimal/best-known (OBK) & \multicolumn{2}{c|}{TL as Table~\ref{table2}}
&  \multicolumn{2}{c|}{TL 1000 seconds} \\
\cline{3-6}
& & Solution & Total time (secs)
& Solution & Total time (secs) \\
\hline
										
4.1	&	429	&	o	&	0.0	&	o	&	0.0	\\
4.2	&	512	&	o	&	0.0	&	o	&	0.0	\\
4.3	&	516	&	o	&	0.0	&	o	&	0.0	\\
4.4	&	494	&	o	&	0.1	&	o	&	0.1	\\
4.5	&	512	&	o	&	0.0	&	o	&	0.0	\\
4.6	&	560	&	o	&	0.2	&	o	&	0.1	\\
4.7	&	430	&	o	&	0.0	&	o	&	0.0	\\
4.8	&	492	&	o	&	0.1	&	o	&	0.1	\\
4.9	&	641	&	o	&	0.1	&	o	&	0.2	\\
4.1	&	514	&	o	&	0.1	&	o	&	0.1	\\
\hline
5.1	&	253	&	o	&	0.1	&	o	&	0.1	\\
5.2	&	302	&	o	&	0.1	&	o	&	0.2	\\
5.3	&	226	&	o	&	0.1	&	o	&	0.0	\\
5.4	&	242	&	o	&	0.1	&	o	&	0.1	\\
5.5	&	211	&	o	&	0.1	&	o	&	0.1	\\
5.6	&	213	&	o	&	0.1	&	o	&	0.0	\\
5.7	&	293	&	o	&	0.1	&	o	&	0.1	\\
5.8	&	288	&	o	&	0.1	&	o	&	0.1	\\
5.9	&	279	&	o	&	0.0	&	o	&	0.1	\\
5.1	&	265	&	o	&	0.0	&	o	&	0.0	\\
\hline
6.1	&	138	&	o	&	0.3	&	o	&	0.2	\\
6.2	&	146	&	o	&	0.2	&	o	&	0.1	\\
6.3	&	145	&	o	&	0.1	&	o	&	0.1	\\
6.4	&	131	&	o	&	0.1	&	o	&	0.1	\\
6.5	&	161	&	o	&	0.4	&	o	&	0.3	\\
\hline
A1	&	253	&	o	&	0.3	&	o	&	0.3	\\
A2	&	252	&	o	&	0.5	&	o	&	0.3	\\
A3	&	232	&	o	&	0.3	&	o	&	0.3	\\
A4	&	234	&	o	&	0.2	&	o	&	0.1	\\
A5	&	236	&	o	&	0.1	&	o	&	0.1	\\
\hline
B1	&	69	&	o	&	0.5	&	o	&	0.4	\\
B2	&	76	&	o	&	0.8	&	o	&	0.8	\\
B3	&	80	&	o	&	0.5	&	o	&	0.5	\\
B4	&	79	&	o	&	1.1	&	o	&	1.2	\\
B5	&	72	&	o	&	0.4	&	o	&	0.5	\\
\hline
C1	&	227	&	o	&	0.4	&	o	&	0.3	\\
C2	&	219	&	o	&	0.7	&	o	&	0.6	\\
C3	&	243	&	o	&	1.2	&	o	&	1.1	\\
C4	&	219	&	o	&	0.2	&	o	&	0.2	\\
C5	&	215	&	o	&	0.4	&	o	&	0.3	\\
D1	&	60	&	o	&	1.1	&	o	&	1.0	\\
D2	&	66	&	o	&	1.6	&	o	&	1.6	\\
D3	&	72	&	o	&	1.6	&	o	&	1.4	\\
D4	&	62	&	o	&	1.9	&	o	&	1.5	\\
D5	&	61	&	o	&	0.5	&	o	&	0.4	\\
\hline
NRE1	&	29	&	o	&	9.2	&	o	&	8.2	\\
NRE2	&	30	&	o	&	43.3	&	o	&	40.0	\\
NRE3	&	27	&	o	&	16.3	&	o	&	17.5	\\
NRE4	&	28	&	o	&	13.7	&	o	&	12.7	\\
NRE5	&	28	&	o	&	8.6	&	o	&	7.8	\\
\hline
NRF1	&	14	&	o	&	15.0	&	o	&	11.3	\\
NRF2	&	15	&	o	&	10.8	&	o	&	8.0	\\
NRF3	&	14	&	o	&	4.2	&	o	&	3.1	\\
NRF4	&	14	&	o	&	20.1	&	o	&	17.8	\\
NRF5	&	13	&	o	&	48.4	&	o	&	41.0	\\
\hline
NRG1	&	176	&	o	&	308.4	&	o	&	1000.0	\\
NRG2	&	154	&	155	&	219.7	&	o	&	335.1	\\
NRG3	&	166	&	167	&	565.4	&	167	&	1000.0	\\
NRG4	&	168	&	o	&	821.4	&	o	&	1000.0	\\
NRG5	&	168	&	o	&	475.4	&	o	&	1000.0	\\
\hline
NRH1	&	63	&	65	&	345.0	&	64	&	1000.0	\\
NRH2	&	63	&	64	&	501.1	&	64	&	1000.0	\\
NRH3	&	59	&	o	&	773.9	&	o	&	1000.0	\\
NRH4	&	58	&	59	&	619.2	&	o	&	1000.0	\\
NRH5	&	55	&	o	&	342.0	&	o	&	1000.0	\\
\hline
\end{tabular}
\caption{Computational results: Cplex experiments}
\label{table4}
\end{table}
\normalsize

\begin{table}[hbtp!]
\centering
\renewcommand{\tabcolsep}{1mm} 
\renewcommand{\arraystretch}{1.0} 
\begin{tabular}{|l|l|c|c|c|}
\hline 
 & Approach & Average & Number of & Average total\\
&  & \% deviation & OBK solutions & time (secs) \\
\hline
& This paper & 0.02 & 64/65 & 94.6 \\
All 65  problems & Cplex: TL as Table~\ref{table2} & 0.12 & 60/65 & 79.7 \\
& Cplex: TL 1000 seconds & 0.06 & 62/65 & 146.4 \\
\hline
& This paper & 0.18 & 8/9 & 527.7 \\
9 NRG and NRH problems & Cplex: TL as Table~\ref{table2} & 0.79 & 5/9 & 528.0 \\
& Cplex: TL 1000 seconds & 0.42 & 6/9 & 1000.0 \\
\hline
\end{tabular}
\caption{Comparison: this paper versus Cplex alone}
\label{table5}
\end{table}

\begin{table}[hbtp!]
\centering
\renewcommand{\tabcolsep}{1mm} 
\renewcommand{\arraystretch}{1.0} 
\begin{tabular}{|l|c|c|}
\hline 
Instance & Optimal solution value 
& Total time (secs) \\
\hline
NRG1	& 176	& 2071.4 \\
NRG3	& 166 &	9350.8  \\
NRG4	& 168	& 4338.8 \\
NRG5	& 168	& 16484.3 \\
\hline
\end{tabular}
\caption{Optimal solution values as obtained using Cplex}
\label{tableopt}
\end{table}

\subsection{Comparison with other heuristics}

In order to compare our results with previous results in the literature we have taken the detailed results given by various authors and computed percentage deviation from the 
optimal/best-known value, OBK, as shown in Table~\ref{table2}. In other words using the solution values given by authors in their papers we computed 100(solution value - OBK)/OBK for each individual problem  and then averaged the percentage deviations. Note here that whilst previous authors may have used this percentage deviation approach their OBK values may be different from those shown in Table~\ref{table2} (since obviously best-known values may be updated over time).

Table~\ref{table3} shows this comparison, where  in each case we give the average percentage deviation, the number of solutions equal 
to the optimal/best-known solution (OBK) and the average computation total time (in seconds), as calculated 
using solution details in the papers cited.
In a similar manner as in Table~\ref{table5} we have split Table~\ref{table3} into two sections, one dealing with all 65 problems, one dealing with the nine larger NRG and NRH problems that were challenging 
for Cplex to solve.

As we would expect different authors have used different hardware, as shown in Table~\ref{table3}. This makes a direct 
comparison between computation times very difficult. However a \emph{\textbf{very rough approximation}} of the relative speed of the different hardware used can be gained by looking at the Mflops/s values as from https://www.cpubenchmark.net/singleCompare.php~\cite{passmark} or given in Dongarra~\cite{dongarra}, although that latter reference has not been updated since 2014. Utilising these sources we give in  Table~\ref{table3}  our \emph{\textbf{best estimate}} of the scaling factor needed to convert the computation times reported on the differing hardware cited  into a normalised time (seconds on our hardware). Here the scaling factor indicates how much faster our hardware is than the hardware associated with other work.
So for example in Table~\ref{table3} the average total time for all 65 problems for Aickelin~\cite{aickelin02}   is 1179.5 seconds  (when using a DECstation 5000/240), which with a scaling factor of 730.5 equates to a normalised time of 1179.5/730.5 =1.6  seconds on the hardware we used.

With regard to Table~\ref{table3}:
\begin{compactitem}
\item Aickelin~\cite{aickelin02} and Beasley and Chu~\cite{beasley96} used ten trials of their algorithms. So for both these papers 
in Table~\ref{table3} we have used the best solution found over the ten trials in calculating percentage deviation, and the time seen is the total time for these ten trials.
 \item For Lan et al~\cite{lan07} the detailed times given in  their paper are the times as to when the final best solution was first found, rather than total times.
 For this reason we used Lan~\cite{lan04}, which does give total times in order to present the values for Lan in Table~\ref{table3}.
 The results given in
Table~\ref{table3} are for the best performing of three different variants of their approach (the variant Meta-RaPS with randomized priority rules).
\item For Reyes and Araya~\cite{reyes21} the times given in their paper appear to be the average time for each trial, where they used 30 trials. So in Table~\ref{table3} we give the total time for 30 trials and in calculating percentage deviation we used the best solution found over these 30 trials.
\item For Soto et al~\cite{soto17} the times given in their paper are the average time for each trial~\cite{private22}, where  they used 30 trials. So in Table~\ref{table3}   we give the total time for 30 trials and in calculating percentage  deviation we used the best solution found over these  30 trials.
\end{compactitem}

\noindent  For Caprara et al~\cite{caprara99} the detailed times given in  their paper are the times as to when the final best solution was first found. As such we have no detailed information as to the total time taken, which is the value given for the other papers cited. So, in order to make an appropriate comparison, we give in Table~\ref{tablec} results for our approach using the time at which the final best solution was first found.

 \emph{\textbf{We should stress here that the papers considered in Table~\ref{table3} and Table~\ref{tablec}
represent, based upon our literature review, the most effective heuristics for the non-unicost  set covering problem given previously in the literature.}} 

Considering Table~\ref{table3} it is clear that, for the 65 test problems examined, the heuristic presented in this paper gives good results in terms of solution quality as compared with the seven other approaches from the literature shown in that table.
Our approach has the second lowest average percentage deviation, and the second highest number of OBK solutions, both over the 65 test problems and over the nine NRG and NRH problems
that Cplex could not solve easily. 
In  terms of solution quality then from
Table~\ref{table3} and Table~\ref{tablec} only 
Lan~\cite{lan04,lan07} and  Caprara et al~\cite{caprara99}  out-perform our approach. This is due to the fact that both algorithms find one better solution over the 65 test problems. 

Obviously, in judging the worth of a heuristic for a specific problem, computation time must be taken into account. Clearly it is a judgement call but it seems reasonable to suggest that, in terms of the non-unicost set covering problem, 
our heuristic is better than six of the eight other heuristics we examined,  slightly worse than that of Lan~\cite{lan04,lan07}, but that all  heuristics for the problem are out-performed by Caprara et al~\cite{caprara99}.

However it is important to note here that, aside from a simple greedy heuristic to find an initial feasible solution, our heuristic made no use of problem-specific knowledge.  This is in sharp contrast to the other heuristics seen in Table~\ref{table3} and Table~\ref{tablec}, which do make use of such knowledge, e.g. in terms of problem-specific preprocessing
to eliminate columns. 

We would stress here that, although the heuristic presented in this paper is not the best-known heuristic for the non-unicost set covering problem, 
we believe that  the work described here illustrates that the potential for using local branching, operating as a stand-alone matheuristic as in the approach described in this paper,
has not been fully exploited in the literature.

\begin{landscape}
\begin{table}[hbtp!]
\footnotesize
\centering
\renewcommand{\tabcolsep}{1mm} 
\renewcommand{\arraystretch}{1.0} 
\begin{tabular}{|c|l|c|c|c|c|c|c|}

\hline 
& Approach & Average & Number of & Average total& Computer/processor & Best estimate &  Average\\
& & \% deviation & OBK solutions & time (secs) & used & scaling factor & normalised time (sec) \\
\hline
& This paper & 0.02 & 64/65 & 94.6 & Intel Core i5-1135G7 2.4GHz & 1.0 & 94.6 \\
& Aickelin~\cite{aickelin02} & 0.13 & 61/65 & 1179.5 & DECstation 5000/240  &730.5 &1.6 \\
All & Beasley and Chu~\cite{beasley96} & 0.07 & 61/65 & 14694.3 & Silicon Graphics Indigo R4000 100MHz & 225.3& 65.2\\
 65 & Lan~\cite{lan04, lan07} & 0 & 65/65 & 1227.1 & Intel Pentium IV 1.7GHz &10.6 & 115.8 \\
problems & Naji-Azimi et al~\cite{naji10} & 0.18 & 53/65 & 118.5 & Intel Core Duo 1.7GHz & 5.1& 23.2 \\
& Reyes and Araya~\cite{reyes21} & 0.25 & 50/65 & 244.7 & Intel Core i7-4700MQ 2.4Ghz & 1.5&163.1 \\
& Soto et al~\cite{soto17} black hole & 1.55 & 35/65 & 184.5 & Intel Quad Core 3.0GHz &  2.1&87.9\\
& Soto et al~\cite{soto17} cuckoo search & 1.01 &35/65 & 151.0 & Intel Quad Core 3.0GHz & 2.1&71.9 \\
\hline
& This paper & 0.18 & 8/9 & 527.7 & Intel Core i5-1135G7 2.4GHz & 1.0 & 527.7 \\
& Aickelin~\cite{aickelin02} & 0.72 & 7/9 & 3830.0 & DECstation 5000/240  &730.5 & 5.2 \\
9 NRG & Beasley and Chu~\cite{beasley96} & 0.35& 7/9 & 36934.4 & Silicon Graphics Indigo R4000 100MHz & 225.3& 163.9 \\
 and & Lan~\cite{lan04, lan07} & 0 & 9/9 & 7764.8 & Intel Pentium IV 1.7GHz &10.6 & 732.5 \\
NRH problems & Naji-Azimi et al~\cite{naji10} & 0.96 & 4/9 & 355.9 & Intel Core Duo 1.7GHz & 5.1& 69.8 \\
& Reyes and Araya~\cite{reyes21} & 1.14 & 3/9 & 1386.0  & Intel Core i7-4700MQ 2.4Ghz & 1.5& 924.0 \\
& Soto et al~\cite{soto17} black hole & 5.26 & 1/9  & 544.5 & Intel Quad Core 3.0GHz &  2.1& 259.3 \\
& Soto et al~\cite{soto17} cuckoo search & 1.78 & 1/9 & 389.4 & Intel Quad Core 3.0GHz & 2.1& 185.4 \\
\hline
\end{tabular}
\caption{Comparison of results}
\label{table3}
\end{table}

\begin{table}[hbtp!]
\footnotesize
\centering
\renewcommand{\tabcolsep}{1mm} 
\renewcommand{\arraystretch}{1.0} 
\begin{tabular}{|c|l|c|c|c|c|c|c|}

\hline 
& Approach & Average & Number of & Average time to & Computer/processor & Best estimate &  Average\\
& & \% deviation & OBK solutions & best solution (secs) & used & scaling factor & normalised time (sec) \\
\hline
All 65& This paper & 0.02 & 64/65 & 39.0 & Intel Core i5-1135G7 2.4GHz & 1.0 & 39.0 \\
problems & Caprara et al~\cite{caprara99} & 0 & 65/65 & 139.9 & DECstation 5000/240 & 730.5 & 0.2 \\
\hline
9 NRG and NRH & This paper & 0.18 & 8/9 & 257.3 & Intel Core i5-1135G7 2.4GHz & 1.0 & 257.3 \\
problems & Caprara et al~\cite{caprara99} & 0 & 9/9 & 670.2 & DECstation 5000/240 & 730.5 & 0.9\\
\hline
\end{tabular}
\caption{Comparison of results: Caprara et al~\cite{caprara99}}
\label{tablec}
\end{table}
\normalsize 
\normalsize 
\end{landscape}


\section{Conclusions}
\label{conc}

In this paper we have presented 
 a heuristic for the non-unicost set covering problem using local branching.
Local branching eliminates the need to define a problem specific search neighbourhood for any particular (zero-one) optimisation problem. It does this by incorporating a generalised Hamming distance neighbourhood into the problem, and this leads naturally to an appropriate neighbourhood search procedure. 

We applied our approach to the non-unicost set covering problem and presented computational results for 65 test problems that have been widely considered in the literature. Our results indicated that our heuristic for the set covering problem using local branching 
is better than six of the eight other heuristics we examined,  slightly worse than that of one heuristic, but that
there is a single heuristic that out-performs all others.

We would  stress here the relative simplicity involved in  creating a local branching based matheuristic 
using the approach given in this paper. Obviously one needs a mathematical formulation of the problem at hand, but aside 
from that the most that is needed is some (relatively simple) problem-specific procedure to generate an initial feasible solution. 

We believe that the work described here illustrates that the potential for using local branching, operating as a stand-alone matheuristic as in the approach described in this paper,
has not been fully exploited in the literature.

It is important to note however that considering the set covering problem alone is insufficient to determine the full potential of
using local branching operating as a stand-alone matheuristic.
Further analysis and evaluation, involving other problems, is
necessary to establish the  effectiveness of the approach proposed in this paper, to see whether it leads to algorithmic success.
 We hope that the work presented here will encourage other researchers to explore using local branching as a stand-alone matheuristic for other problems.

\FloatBarrier
 \clearpage
\newpage
 \pagestyle{empty}
\linespread{1}
\small \normalsize 

\section*{Acknowledgments}
\noindent The author would like to acknowledge the comments made on an earlier version of this paper by anonymous reviewers.

\section*{Conflict of interest statement}
The author has no relevant financial or non-financial interests to disclose.

\bibliography{paper}

\begin{thebibliography}{10}

\bibitem{aickelin02}
U.~Aickelin.
\newblock An indirect genetic algorithm for set covering problems.
\newblock {\em Journal of the Operational Research Society}, 53(10):1118--1126,
  2002.

\bibitem{aranha22}
C.~Aranha, C.~L.~C. Villalon, F.~Campelo, M.~Dorigo, R.~Ruiz, M.~Sevaux,
  K.~Sorensen, and T.~Stutzle.
\newblock Metaphor-based metaheuristics, a call for action: the elephant in the
  room.
\newblock {\em Swarm Intelligence}, 16(1):1--6, 2022.

\bibitem{balas80}
E.~Balas and A.~Ho.
\newblock Set covering algorithms using cutting planes, heuristics, and
  subgradient optimization - a computational study.
\newblock {\em Mathematical Programming Study}, 12(APR):37--60, 1980.

\bibitem{beasley87}
J.~E. Beasley.
\newblock An algorithm for set covering problem.
\newblock {\em European Journal of Operational Research}, 31(1):85--93, 1987.

\bibitem{beasley90}
J.~E. Beasley.
\newblock A lagrangian heuristic for set-covering problems.
\newblock {\em Naval Research Logistics}, 37(1):151--164, 1990.

\bibitem{beasley1990}
J.~E. Beasley.
\newblock {OR-Library: distributing test problems by electronic mail}.
\newblock {\em {Journal of the Operational Research Society}},
  {41}({11}):{1069--1072}, {1990}.

\bibitem{beasley96}
J.~E. Beasley and P.~C. Chu.
\newblock A genetic algorithm for the set covering problem.
\newblock {\em European Journal of Operational Research}, 94(2):392--404, 1996.

\bibitem{birbil03}
S.~I. Birbil and S.~C. Fang.
\newblock An electromagnetism-like mechanism for global optimization.
\newblock {\em Journal of Global Optimization}, 25(3):263--282, 2003.

\bibitem{boschetti23}
M.~A. Boschetti, A.~N. Letchford, and V.~Maniezzo.
\newblock Matheuristics: survey and synthesis.
\newblock {\em International Transactions in Operational Research}, 30(6,
  SI):2840--2866, 2023.

\bibitem{boschetti2009}
M.~A. Boschetti, V.~Maniezzo, M.~Roffilli, and A.~B. Rohler.
\newblock {\em {{Matheuristics: Optimization, simulation and control}}}, volume
  {5818} of {\em {Lecture Notes in Computer Science}}, pages {171--177}.
\newblock {Springer, Berlin}, {2009}.

\bibitem{caprara99}
A.~Caprara, M.~Fischetti, and P.~Toth.
\newblock A heuristic method for the set covering problem.
\newblock {\em Operations Research}, 47(5):730--743, 1999.

\bibitem{caprara00}
A.~Caprara, P.~Toth, and M.~Fischetti.
\newblock Algorithms for the set covering problem.
\newblock {\em Annals of Operations Research}, 98:353--371, 2000.

\bibitem{cplex1210}
{CPLEX Optimizer, version 12.10}.
\newblock {IBM. Available from
  https://www.ibm.com/products/ilog-cplex-optimization-studio last accessed
  April 2023}, {2019}.

\bibitem{dongarra}
J.~J. Dongarra.
\newblock {Performance of various computers using standard linear equations
  software}.
\newblock {Available from http://www.netlib.org/benchmark/performance.ps last
  accessed November 2023}, {2014}.

\bibitem{feo95}
T.~A. Feo and M.~G.~C. Resende.
\newblock Greedy randomized adaptive search procedures.
\newblock {\em Journal of Global Optimization}, 6(2):109--133, 1995.

\bibitem{festa02}
P.~Festa and M.~G.~C. Resende.
\newblock {GRASP: An annotated bibliography}.
\newblock In C.~C. Ribeiro and P.~Hansen, editors, {\em Essays and Surveys in
  Metaheuristics}, volume~15 of {\em Operations Research Computer Science
  Interfaces}, pages 325--367. Kluwer Academic Publishers, 2002.

\bibitem{fischetti03}
M.~Fischetti and A.~Lodi.
\newblock Local branching.
\newblock {\em Mathematical Programming}, 98(1-3):23--47, 2003.

\bibitem{glover90}
F.~Glover.
\newblock Tabu search - a tutorial.
\newblock {\em Interfaces}, 20(4):74--94, 1990.

\bibitem{hansen01}
P.~Hansen and N.~Mladenovic.
\newblock Variable neighborhood search: principles and applications.
\newblock {\em European Journal of Operational Research}, 130(3):449--467,
  2001.

\bibitem{kirkpatrick83}
S.~Kirkpatrick, C.~D. Gelatt, and M.~P. Vecchi.
\newblock Optimization by simulated annealing.
\newblock {\em Science}, 220(4598):671--680, 1983.

\bibitem{kumar15}
S.~Kumar, D.~Datta, and S.~K. Singh.
\newblock Black hole algorithm and its applications.
\newblock In A.~T. Azar and S.~Vaidyanathan, editors, {\em Computational
  Intelligence Applications in Modeling and Control}, pages 147--170. Springer
  International Publishing, 2015.

\bibitem{lan04}
G.~Lan.
\newblock {An effective and simple heuristic based on Meta-RaPS to solve
  large-scale set covering problems}.
\newblock {M.S. Thesis, University of Louisville. A electronic copy of this
  thesis can be obtained from the University of Louisville library
  https://library.louisville.edu}, 2004.

\bibitem{lan07}
G.~Lan, G.~W. DePuy, and G.~E. Whitehouse.
\newblock An effective and simple heuristic for the set covering problem.
\newblock {\em European Journal of Operational Research}, 176(3):1387--1403,
  2007.

\bibitem{maniezzo21}
V.~Maniezzo, M.~A. Boschetti, and T.~Stuzle.
\newblock {\em {{Matheuristics: Algorithms and implementations}}}.
\newblock {EURO Advanced Tutorials on Operational Research}. {Springer,
  Berlin}, {2021}.

\bibitem{mladenovic97}
N.~Mladenovic and P.~Hansen.
\newblock Variable neighborhood search.
\newblock {\em Computers \& Operations Research}, 24(11):1097--1100, 1997.

\bibitem{naji10}
Z.~Naji-Azimi, P.~Toth, and L.~Galli.
\newblock An electromagnetism metaheuristic for the unicost set covering
  problem.
\newblock {\em European Journal of Operational Research}, 205(2):290--300,
  2010.

\bibitem{private22}
R.~Olivares.
\newblock Private communication, 2022.

\bibitem{passmark}
{Passmark Software: CPU Benchmarks}.
\newblock {https://www.cpubenchmark.net/singleCompare.php last accessed
  November 2023}, {2023}.

\bibitem{reyes21}
V.~Reyes and I.~Araya.
\newblock {A GRASP-based scheme for the set covering problem}.
\newblock {\em Operational Research}, 21(4):2391--2408, 2021.

\bibitem{soto17}
R.~Soto, B.~Crawford, R.~Olivares, J.~Barraza, I.~Figueroa, F.~Johnson,
  F.~Paredes, and E.~Olguin.
\newblock Solving the non-unicost set covering problem by using cuckoo search
  and black hole optimization.
\newblock {\em Natural Computing}, 16(2):213--229, 2017.

\bibitem{yang09}
X.-S. Yang and S.~Deb.
\newblock {Cuckoo search via L\'evy flights}.
\newblock In A.~Abraham, F.~Herrera, A.~Carvalho, and V.~Pai, editors, {\em
  2009 World Congress on Nature \& Biologically Inspired Computing (NABIC
  2009)}, pages 210--214, 2009.

\end{thebibliography}

\end{document}